\documentclass[11pt, reqno]{amsart}
\usepackage{amsmath}
\usepackage{amssymb}
\usepackage{xcolor}
\usepackage[colorlinks, linkcolor=purple, citecolor=blue, urlcolor=blue, hypertexnames=false]{hyperref}
\usepackage{geometry}
\usepackage{graphicx}

\numberwithin{equation}{section}

\newtheorem{thm}{Theorem}[section]
\newtheorem{prop}{Proposition}[section]

\newtheorem{lem}{Lemma}[section]
\newtheorem{cor}{Corollary}[section]

\newtheorem{rem}{Remark}[section]

\author[S. Kim, T. Saanouni and I. Seo]{Seongyeon Kim, Tarek Saanouni and Ihyeok Seo}

\address[S. Kim]{Department of Mathematics Education, Jeonju University, Jeonju 55069, Republic of Korea}
\email{\it \color{blue}{sy\_kim@jj.ac.kr}}

\address[T. Saanouni]{Departement of Mathematics, College of Science, Qassim University, Buraydah, Kingdom of Saudi Arabia}
\email{\it\color{blue}{t.saanouni@qu.edu.sa}}

\address[I. Seo]{Department of Mathematics, Sungkyunkwan University, Suwon 16419, Republic of Korea}
\email{\it \color{blue}{ihseo@skku.edu}}


\title[The energy-critical inhomogeneous Hartree equation]{Remarks on the well-posedness of the energy-critical inhomogeneous Hartree equation}

\begin{document}

\begin{abstract}
We study the energy-critical inhomogeneous Hartree equation in space dimensions three and higher. Previous local well-posedness results left open the parameter regime where the inhomogeneity exponent is small and the Riesz potential exponent is either small or large. We establish local well-posedness for a new range of parameters, thereby substantially filling the remaining open parameter regime. In particular, our result completely resolves the remaining gap in dimensions $5$ and $6$. 
\end{abstract}

\maketitle

\section{Introduction}\label{sec1}
In this paper, we consider the Cauchy problem for the inhomogeneous Hartree equation
\begin{equation}\label{S}
\begin{cases}
i\partial_t u + \Delta u = \epsilon (I_\alpha * |x|^{-b} |u|^p) |x|^{-b} |u|^{p-2} u, \quad (x,t) \in \mathbb{R}^n \times \mathbb{R}, \\
u(x,0) = u_0(x),
\end{cases}
\end{equation}
where $p\ge2$, $b >0$, and $\epsilon = \pm 1$. The case $\epsilon = 1$ is referred to as the \textit{defocusing case}, while $\epsilon = -1$ corresponds to the \textit{focusing case}.
The Riesz potential $I_\alpha$ on $\mathbb{R}^n$ is defined  by
$$ I_{\alpha} (x):= \frac{\Gamma\left(\frac{n-\alpha}{2}\right)}{\Gamma\left(\frac{\alpha}{2}\right) \pi^{n/2} 2^\alpha} |x|^{-(n-\alpha)}, \quad 0 < \alpha < n.$$
This equation appears in various physical models, where the factor $|x|^{-b}$ describes an inhomogeneity of the underlying medium (see \cite{Gi}). Particular instances of this equation arise, for example, in the mean-field description of large systems of non-relativistic atoms and molecules, as well as in the propagation of electromagnetic waves in plasmas. We refer to \cite{BC,fl,r} for more details.

We first recall the critical Sobolev exponent associated with \eqref{S}. Equation \eqref{S} is invariant under the scaling
$$ u(x,t) \mapsto u_\delta(x,t) := \delta^{\frac{2 - 2b + \alpha}{2(p-1)}} u(\delta x, \delta^2 t), \quad \delta > 0.$$
Under this scaling, the homogeneous $\dot H^1$ norm of the initial data scales as
$$\|u_\delta(\cdot, 0)\|_{\dot{H}^1} = \delta^{1 - \frac{n}{2} + \frac{2 - 2b + \alpha}{2(p - 1)}} \|u_0\|_{\dot{H}^1}.$$
Hence the $\dot H^1$ norm is invariant precisely when
\begin{equation}\label{power}
p = 2+\frac{\alpha-n+4-2b}{n - 2}.
\end{equation}
We refer to this exponent as the energy-critical power and call \eqref{S}  the \emph{energy-critical inhomogeneous Hartree equation} in this case.

\begin{figure}
\centering
\includegraphics[width=1.0\linewidth]{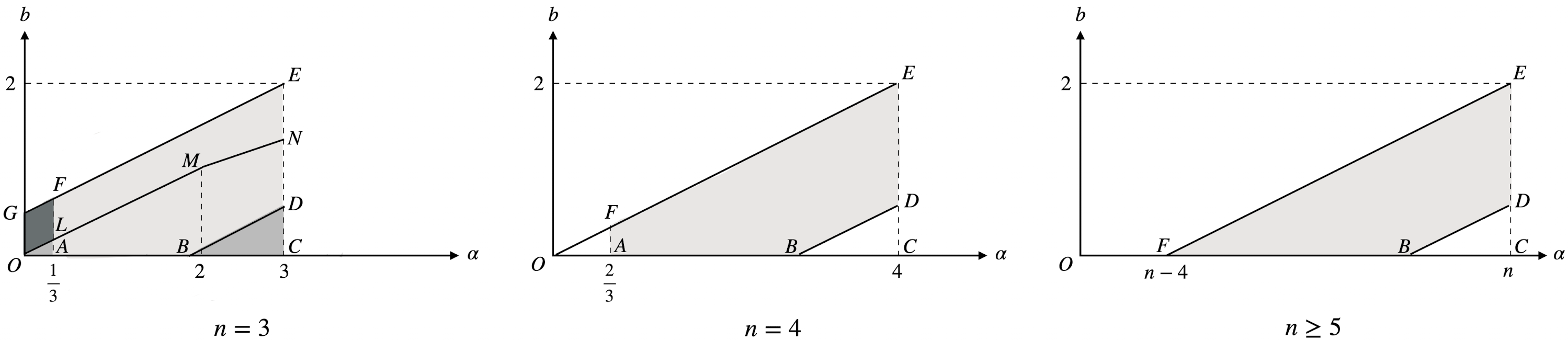}
\caption{Known regions of $(\alpha, b)$ for the local well-posedness of the energy-critical \eqref{S}. Here, $B=( \frac{n-4+\sqrt{9n^2-8n+16}}4,0)$, $D=(n,\frac{3n+4-\sqrt{9n^2-8n+16}}{8})$, $F=(\max\{\frac{n-2}{3},n-4\}, \max\{0, \frac{5-n}{3}\})$, $G=(0,\frac{1}{2})$, $L(\frac13,\frac16)$, $M(2,1)$, $N(3,\frac43)$.} \label{figure}
\end{figure}

For $n\ge3$, since $p\ge2$, it follows from \eqref{power} that the admissible range of the parameters $(\alpha,b)$ in the energy-critical regime is given by
\begin{equation}\label{range}
\max\{0,n-4\}<\alpha<n,\quad 0< b\le\frac{\alpha-n+4}{2}.
\end{equation}
The line $b=(\alpha-n+4)/2$ passes through the points $E$ and $F$ 
in Figure \ref{figure}.
Therefore, the range corresponds to the region lying below this line in the figure.
Within the full parameter range \eqref{range}, energy-critical local well-posedness remains unsolved for $n\ge4$ in the regime where $b$ is small and $\alpha$ is either small or large. In this paper, we address this remaining regime for all dimensions $n\ge4$. In particular, for $n=5,6$, our result completely fills the remaining gap.

To clarify this point, we first summarize the parameter ranges addressed in the previous results. 
When $b= 0$ and $p = 2$, 
the equation \eqref{S} reduces to the classical Hartree (or Choquard) equation, in which case 
\eqref{power} implies $\alpha = n-4$. Local well-posedness in this case was established 
by J. Ginibre and G. Velo \cite{gvl} in 1980.
In the more general case $p>2$, the parameter range for $\alpha$ is given by $\max\{0, n-4\}<\alpha<n$, and this case was subsequently resolved in \cite{AR}. Consequently, the homogeneous case $b=0$ has been completely settled.

When $b>0$, the inhomogeneous factor $|x|^{-b}$ in the nonlinearity introduces analytical difficulties in the critical regime. To overcome these difficulties, S. Kim \cite{sk} employed Sobolev–Lorentz spaces to absorb the singularity and used the corresponding Strichartz-type estimates. 
The parameter range obtained in \cite[Theorem 1.2]{sk} is given by 
\begin{equation}\label{kim}
\max\Big\{\frac{n-2}{3}, n-4\Big\} <\alpha < n,\, \max\big\{0, \frac{\alpha}2-\frac{n-4+\sqrt{9n^2-8n+16}}8\big\}<b \le \frac{\alpha-n+4}{2}.
\end{equation}
The line $b=\frac{\alpha}2-\frac{n-4+\sqrt{9n^2-8n+16}}8$ passes through the points $B$ and $D$ in Figure \ref{figure}.
Hence, the parameter range \eqref{kim} corresponds to the lightest shaded region in Figure \ref{figure}.
Therefore, the region $BCD$ remains open, together with the region $OAFG$ when $n=3$ and the region $OAF$ when $n=4$.

Subsequently, the remaining open region in the three-dimensional case was almost completely resolved in the works \cite{GX,SP}, except for the borderline segment $(G,F]$ in Figure \ref{figure}. More precisely, the parameter range obtained in \cite[Theorem 1.1]{GX} is given by 
\begin{equation}\label{gx}
0<\alpha<3,\quad 0 < b \leq \min\Big\{\frac{\alpha+1}{3}, \frac{\alpha}{2}\Big\}. 
\end{equation}
In Figure \ref{figure}, the lines $b=(\alpha+1)/3$ and $b=\alpha/2$ correspond to the segments $MN$ and $OM$, respectively. Hence, the parameter range \eqref{gx} corresponds to the region lying below these segments, and the remaining open region is reduced to $OLFG$, corresponding to
\[
0<\alpha\leq\frac13,
\qquad
\frac{\alpha}{2}<b\leq\frac{\alpha+1}{2}.
\]
This remaining region was subsequently covered in \cite{SP}, except for the borderline segment $(G,F]$, where local well-posedness was established for the range
\begin{equation}\label{sp}
0<\alpha<1,
\qquad
\frac{\alpha}{2} < b < \frac{\alpha+1}{2}.
\end{equation}

However, for $n\ge4$, there had been no further progress since the work \cite{sk}. We aim to address the remaining open regions in these higher-dimensional cases. The theorem below substantially fills the remaining gaps and, in particular, completely resolves the remaining gap when $n=5,6$. For a more detailed discussion, see the paragraph below Corollary \ref{cor}.

\begin{thm}\label{thm}
Let $n\ge3$ and $u_0 \in H^1(\mathbb{R}^n)$.
Assume that
\begin{equation}\label{as}
\max\{0,n-4\}<\alpha<n \quad \text{and} \quad   
0<b \le \frac{\alpha-n+4}{n}.
\end{equation}
Then there exist $T>0$ and a unique solution
$$u\in C([0,T]; H^1(\mathbb R^n)) \cap L^q([0,T];H^{1,r}(\mathbb R^n)),$$
where $(q,r)$ is a Schr\"odinger admissible pair satisfying
$$2\le q\le \infty, \qquad \frac{2}{q}+\frac{n}{r}=\frac{n}{2},$$
to the inhomogeneous Hartree equation \eqref{S} with the energy-critical power $p = 2+\frac{\alpha-n+4-2b}{n-2}$.
Moreover, the solution depends continuously on the initial data.
\end{thm}

\begin{rem}
We remark that, in the three-dimensional case, our result improves the upper bound for $b$ in condition \eqref{gx} to $(\alpha+1)/3$, while lowering the lower bound for $b$ in condition \eqref{sp} to $0$.
Throughout the proof, one can also easily see that Theorem \ref{thm} remains valid in the homogeneous case $b=0$. In this case, the energy-critical condition reduces to $\alpha = n-4$ with $0<\alpha<n$, and one recovers the known local well-posedness results for the classical energy-critical Hartree equation.
\end{rem}

\begin{cor}\label{cor}
Under the assumptions of Theorem \ref{thm}, if $\|u_0\|_{H^1}$ is sufficiently small, then the corresponding solution is global and scatters in $H^1$, i.e., there exists $\phi\in H^1$ such that
$$\lim_{t\to\infty} \|u(t)-e^{-it\Delta}\phi\|_{H^1}=0.$$
\end{cor}

\begin{figure}
\centering
\includegraphics[width=0.95\linewidth]{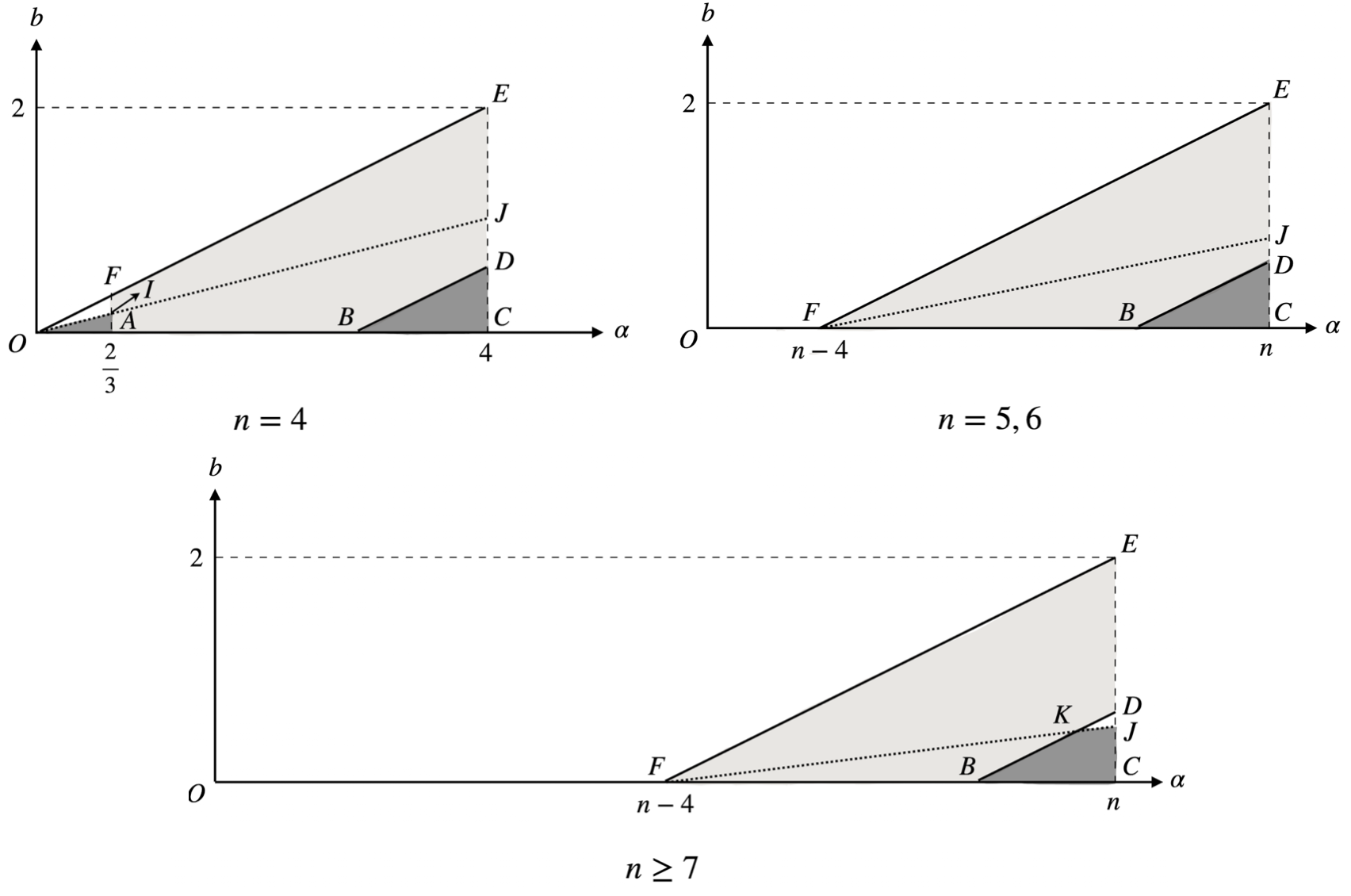}
\caption{Regions of $(\alpha, b)$ for Theorem \ref{thm} when $n\ge4$. 
Here, $B=( \frac{n-4+\sqrt{9n^2-8n+16}}4,0)$, $D=(n,\frac{3n+4-\sqrt{9n^2-8n+16}}{8})$, $I=(\frac23,\frac16)$, $J=(n, \frac4n)$, $F=(\max\{\frac{n-2}{3},n-4\}, \max\{0, \frac{5-n}{3}\})$.} \label{figure2}
\end{figure}

In Figure \ref{figure2}, the line $b=\frac{\alpha-n+4}{n}$ in \eqref{as} corresponds to  the segment $OJ$ when $n=4$, and the segment $FJ$ when $n\geq5$. Therefore, the range \eqref{as} corresponds to the region lying below these segments in the figure, and newly fills the darker shaded region in Figure \ref{figure2}. As a consequence, the only remaining open regions are $OIF$ when $n=4$, and $KJD$ when $n\geq7$. In particular, when $n=5,6$, the remaining region $BCD$ is completely resolved.

The previous approaches \cite{sk,SP} to the critical regime require working in refined function spaces rather than the usual Lebesgue spaces, together with the corresponding Strichartz-type estimates, thereby controlling the singular inhomogeneous factor at the cost of increased technical complexity. In contrast, we observe that by exploiting a weighted Sobolev embedding of Caffarelli-Kohn-Nirenberg type (Lemma \ref{ckn}) to handle the singular factor, one can treat the critical case using only classical Strichartz estimates (Lemma \ref{str}). We also note that this approach is more flexible and extendable than the one in \cite{GX}, which relies on the Hardy inequality.

The rest of this paper is organized as follows.
Section \ref{sec2} is devoted to the derivation of the nonlinear estimates, obtained by combining the Hardy-Littlewood-Sobolev inequality with a weighted Sobolev embedding of Caffarelli-Kohn-Nirenberg type. In Section \ref{sec3}, we use these nonlinear estimates together with the classical Strichartz estimates to prove Theorem \ref{thm}.

Throughout the paper, $C$ denotes a positive constant that may vary from line to line. We write $A \lesssim B$ to mean $A \leq C B$ for some unspecified constant $C > 0$.

\section{Nonlinear estimates}\label{sec2}
In this section, we derive estimates for the nonlinearity appearing in \eqref{S}, which will play a crucial role in establishing the existence of solutions in the next section.
To state these estimates, we first introduce some notation.

We define the nonlinearity
$$\mathcal N[u]= |x|^{-b}|u|^{p-2}(I_\alpha * |\cdot|^{-b}|u|^{p})u$$
for $0<\alpha<n$, $b>0$, $p>2$.
We also define the set of Schr\"odinger admissible pairs
$$\mathcal A = \Big\{(q,r):\, 2 \le q \le \infty, \,\,\, \frac2q + \frac{n}{r}= \frac{n}2\Big\},$$
and the associated norms
$$\|u\|_{S(I)}:= \sup_{(q,r)\in\mathcal A} \|u\|_{L_t^q(I;L_x^r)}$$
and 
$$\|u\|_{S'(I)}:= \inf_{(\tilde q,\tilde r)\in\mathcal A} \|u\|_{L_t^{\tilde q'}(I;L_x^{\tilde r'})}.$$

With the above notation, we obtain the following estimates for the nonlinearity.

\begin{prop}\label{non}
Let \(n\ge3\) and  $0<\alpha<n$.
Assume that
\begin{equation}\label{ass}
0< b\le\frac{\alpha-n+4}{n}.
\end{equation}
Then there exist $(q,r) \in \mathcal{A}$ and $(\tilde q, \tilde r) \in \mathcal{A}$ such that
\begin{equation}\label{non1}
\big\|\nabla\mathcal{N}[u]\big\|_{S'(I)} \lesssim \big\|\nabla u \big\|_{S(I)}^{2p-1}
\end{equation}
and
\begin{equation}\label{non2}
\big\|\mathcal{N}[u]-\mathcal{N}[v]\big\|_{S'(I)} \lesssim \Big(\big\|\nabla u\big\|_{S(I)}^{2p-2} +\big\|\nabla v\big\|_{S(I)}^{2p-2}\Big)\|u-v\|_{S(I)}
\end{equation}
with $p = 2+\frac{\alpha-n+4-2b}{n-2}$.
\end{prop}

\subsection*{Proof}
It suffices to show that there exist pairs $(q,r), (\tilde q , \tilde r )\in \mathcal{A}$ such that the estimates \eqref{non1} and \eqref{non2} hold for these particular pairs, namely
\begin{equation}\label{non11}
\|\nabla\mathcal{N}[u]\|_{L_t^{\tilde q'}(I;L_x^{\tilde r'})} \leq C \|\nabla u\|_{L_t^{q}(I;L_x^r)}^{2p-1}
\end{equation}
and
\begin{equation}\label{non22}
\|\mathcal{N}[u]-\mathcal{N}[v]\|_{L_t^{\tilde q'}(I;L_x^{\tilde r'})} \leq C \|\nabla u\|_{L_t^q(I;L_x^r)}^{2p-2} \|u-v\|_{L_t^q(I;L_x^r)},
\end{equation}
under the assumption \eqref{ass}. 

First, let $(q,r), (\tilde q , \tilde r )\in \mathcal{A}$, namely 
\begin{equation}\label{ad1}
0 \le \frac1q \le \frac12 ,\quad \frac2q + \frac{n}r = \frac{n}2 ,\quad
0 \le \frac1{\tilde q} \le \frac12 ,\quad \frac2{\tilde q} + \frac{n}{\tilde r} = \frac{n}2.
\end{equation}
We also set 
\begin{equation}\label{q}
\frac{1}{\tilde q'} = \frac{2p-1}q.
\end{equation}

We now determine the conditions on the pairs under which \eqref{non11} holds, using \eqref{q}.
Taking the intersection of these conditions with the Schr\"odinger admissibility condition \eqref{ad1}, we obtain \eqref{ass}.
Next, we show that \eqref{non22} follows from a part of the argument used to derive \eqref{non11} and does not impose any additional condition.
Consequently, \eqref{non22} also holds under \eqref{ass}.

To find the conditions for \eqref{non11}, we first note that 
\begin{align*}
&\|\nabla\mathcal{N}[u]\|_{L_t^{\tilde q'}(I;L_x^{\tilde r'})}\\
&\lesssim \big\||x|^{-b-1}|u|^{p-1}(I_\alpha\ast|\cdot|^{-b}|u|^p)\big\|_{L_t^{\tilde q'}(I;L_x^{\tilde r'})}+\big\||x|^{-b}|u|^{p-1}(I_\alpha \ast|\cdot|^{-b-1}|u|^{p})\big\|_{L_t^{\tilde q'}(I;L_x^{\tilde r'})}\\
&\,+ \big\| |x|^{-b}|u|^{p-2}|\nabla u|(I_\alpha\ast|\cdot|^{-b}|u|^p) \big\|_{L_t^{\tilde q'}(I;L_x^{\tilde r'})}+ \big\| |x|^{-b}|u|^{p-1}(I_\alpha\ast|\cdot|^{-b}|u|^{p-1}|\nabla u|) \big\|_{L_t^{\tilde q'}(I;L_x^{\tilde r'})}\\
& =: A_1+A_2+A_3+A_4.
\end{align*}
To bound the terms $A_1$–$A_4$ by $\|\nabla u\|_{L_t^{q}(I;L_x^r)}^{2p-1}$,
we make use of the following two inequalities.
First, we recall the Hardy–Littlewood–Sobolev inequality (\cite{el} and \cite[Corollary 2.14]{tsa}), stated in Lemma \ref{hls}, which is used to control the nonlocal convolution term.
Second, we use a weighted Sobolev embedding of Caffarelli–Kohn–Nirenberg type (\cite{sgw,csl}), given in Lemma \ref{ckn}, to handle the inhomogeneous factor $|x|^{-b}$.

\begin{lem}\label{hls}
Let $n\ge1$, $0 <\alpha < n$ and $1<q,r,s<\infty$.
If $\frac1r +\frac1s =\frac1q+\frac{\alpha}{n}$, then
$$\|f(I_\alpha*g)\|_{L^{q}}\leq C \|f\|_{L^r}\|g\|_{L^s}.$$
\end{lem}

\begin{lem}\label{ckn}
Let $n\geq1$. If
$$1< p\leq q<\infty, \quad -\frac{n}{q}<b\leq 0 \quad \text{and} \quad -b-1=\frac{n}{q}-\frac{n}{p},$$ 
then
$$\||x|^{b}f\|_{L^q}\leq C\|\nabla  f\|_{L^p}.$$
\end{lem}

We begin by estimating the first term $A_1$.  Applying Lemma \ref{hls}, H\"older's inequality in time together with \eqref{q},
and Lemma \ref{ckn}, we obtain
\begin{align*}
A_1 &= \big\||x|^{-b-1}|u|^{p-1} (I_\alpha\ast|\cdot|^{-b}|u|^p)\big\|_{L_t^{\tilde q'}(I;L_x^{\tilde r'})} \\
&\qquad \qquad \qquad \quad  \lesssim \big\||x|^{-b-1}|u|^{p-1}\big\|_{L_t^{\frac{q}{p-1}}(I;L_x^{r_1})} \big\||x|^{-b}|u|^p\big\|_{L_t^{\frac{q}{p}}(I;L_x^{r_2})}\\
&\qquad \qquad \qquad \quad =\big\| |x|^{-\frac{b+1}{p-1}}u \big\|^{p-1}_{L_t^{q}(I;L_x^{(p-1)r_1})}\big\| |x|^{-\frac{b}{p}}u \big\|^{p}_{L_t^{q}(I;L_x^{pr_2})}\lesssim \big\| \nabla u \big\|_{L_t^{q}(I;L_x^{r})}^{2p-1},
\end{align*} 
provided that 
\begin{equation}\label{c7}
0<\frac1{r_1},\,\frac1{r_2}<1, \quad  \frac1{r_1}+\frac1{r_2}=\frac{1}{\tilde r'}+\frac{\alpha}{n},
\end{equation}
\begin{equation}\label{c8}
0<\frac{1}{(p-1)r_1}\leq \frac1r<1,\quad  0\leq\frac{b+1}{p-1}<\frac{n}{(p-1)r_1} ,\quad  \frac{b-p+2}{p-1}=\frac{n}{(p-1)r_1}-\frac{n}{r},
\end{equation}
\begin{equation}\label{c9}
0<\frac{1}{pr_2}\leq \frac1r < 1,\quad 0 \leq \frac{b}{p}<\frac{n}{pr_2},\quad \frac{b-p}{p}=\frac{n}{pr_2}-\frac{n}{r}.
\end{equation}

We next turn to the second term $A_2$. This term can be treated in essentially the same way as $A_1$. Indeed,
\begin{align*}
A_2&=\big\| |x|^{-b}|u|^{p-1}(I_\alpha \ast|\cdot|^{-b-1}|u|^{p})\big\|_{L_t^{\tilde q'}(I;L_x^{\tilde r'})} \\
&\qquad \qquad \quad \quad \lesssim \big\| |x|^{-b}|u|^{p-1}\|_{L_t^{\frac{q}{p-1}}(I;L_x^{r_3})} \big\| |x|^{-b-1}|u|^p \big\|_{L_t^{\frac{q}{p}}(I;L_x^{r_4})}\\
&\qquad \qquad \quad \quad = \big\| |x|^{-\frac{b}{p-1}}u\|^{p-1}_{L_t^{q}(I;L_x^{(p-1)r_3})} \big\| |x|^{-\frac{b+1}{p}}u \big\|^p_{L_t^{q}(I;L_x^{pr_4})}\lesssim \big\| \nabla u \big\|^{2p-1}_{L_t^{q}(I;L_x^{r})},
\end{align*}
under the conditions
\begin{equation}\label{c10}
0<\frac1{r_3},\, \frac{1}{r_4}<1, \quad \frac{1}{r_3}+\frac{1}{r_4}=\frac1{\tilde r'}+\frac{\alpha}{n},
\end{equation}
\begin{equation}\label{c11}
0<\frac1{(p-1)r_3}\leq \frac1r < 1, \quad 0\leq \frac{b}{p-1}<\frac{n}{(p-1)r_3}, \quad \frac{b-p+1}{p-1}=\frac{n}{(p-1)r_3}-\frac{n}{r},
\end{equation}
\begin{equation}\label{c12}
0<\frac{1}{pr_4}\leq \frac1r < 1, \quad 0\leq \frac{b+1}{p}<\frac{n}{pr_4}, \quad \frac{b-p+1}{p}=\frac{n}{pr_4}-\frac{n}{r}.
\end{equation}

We now estimate the third term $A_3$. Applying Lemma \ref{hls}, H\"older's inequality in time together with \eqref{q}, and Lemma \ref{ckn},  we obtain
\begin{align*}
A_3 &=\big\| |x|^{-b}|u|^{p-2}\,|\nabla u| (I_\alpha \ast|\cdot|^{-b}|u|^{p}) \big\|_{L_t^{\tilde q'}(I;L_x^{\tilde r'})}\\
&\quad \quad \quad \quad \lesssim \big\| |x|^{-b}|u|^{p-2}\, \nabla u \big\|_{L_t^{\frac{q}{p-1}}(I;L_x^{r_1})} \big\||x|^{-b}|u|^p \big\|_{L_t^{\frac{q}{p}}(I;L_x^{r_2})}\\
&\quad \quad \quad \quad \leq \big\| |x|^{-b}|u|^{p-2} \big\|_{L_t^{\frac{q}{p-2}}(I;L_x^{r_5})} \big\| |x|^{-b}|u|^p \big\|_{L_t^{\frac{q}{p}}(I;L_x^{r_2})} \| \nabla u \|_{L_t^{q}(I;L_x^{r})}\\
&\quad \quad \quad \quad  =\big\| |x|^{-\frac{b}{p-2}}u \big\|^{p-2}_{L_t^{q}(I;L_x^{(p-2)r_5})} \big\| |x|^{-\frac{b}{p}}u \big\|^p_{L_t^{q}(I;L_x^{pr_2})}\| \nabla u \|_{L_t^{q}(I;L_x^{r})}\lesssim \big\| \nabla u \big\|_{L_t^{q}(I;L_x^{r})}^{2p-1},
\end{align*}
under the conditions \eqref{c7}, $1/{r_1}=1/{r_5}+1/{r}$, \eqref{c9}, and
\begin{equation}\label{c2}
0<\frac{1}{(p-2)r_5} \leq \frac1r< 1,\quad 
0\leq \frac{b}{p-2}<\frac{n}{(p-2)r_5}, \quad \frac{b-p+2}{p-2}=\frac{n}{(p-2)r_5}-\frac{n}{r}.
\end{equation}
Note that the condition $1/{r_1}=1/{r_5}+1/{r}$ follows from combining the last conditions in \eqref{c8} and \eqref{c2}, and hence is redundant.

Finally, we estimate the last term $A_4$, 
which can be handled in the same way as $A_3$. This yields 
\begin{align*}
A_4 &=\big\| |x|^{-b}|u|^{p-1}(I_\alpha\ast|\cdot|^{-b}|u|^{p-1}|\nabla u|) \big\|_{L_t^{\tilde q'}(I;L_x^{\tilde r'})} \\
&\quad \quad \,\, \lesssim \big\| |x|^{-b}|u|^{p-1} \big\|_{L_t^{\frac{q}{p-1}}(I;L_x^{r_3})} \big\| |x|^{-b}|u|^{p-1}\nabla u \big\|_{L_t^{\frac{q}{p}}(I;L_x^{r_4})}\\
&\quad \quad \,\, \leq \big\| |x|^{-\frac{b}{p-1}} u \big\|^{p-1}_{L_t^{q}(I;L_x^{(p-1)r_3})} \big\| |x|^{-\frac{b}{p-1}}u \big\|^{p-1}_{L_t^{q}(I;L_x^{(p-1)r_6})} \| \nabla u \|_{L_t^{q}(I;L_x^{r})} \lesssim \big\| \nabla u \big\|_{L_t^{q}(I;L_x^{r})}^{2p-1},
\end{align*}
under the conditions \eqref{c10}, $1/{r_4}=1/{r_6}+1/r$, \eqref{c11}, and,
\begin{equation}\label{c6}
0<\frac1{(p-1)r_6}\leq \frac1r< 1, \quad 0\leq \frac{b}{p-1}<\frac{n}{(p-1)r_6}, \quad \frac{b-p+1}{p-1}=\frac{n}{(p-1)r_6}-\frac{n}{r}.
\end{equation}
As before, the condition $1/{r_4}=1/{r_6}+1/r$ is redundant, since it follows from the last conditions in  \eqref{c12} and \eqref{c6}.

We now eliminate the auxiliary indices $r_1$, $\cdots$, $r_6$ from the conditions \eqref{c7}–\eqref{c6} required for \eqref{non11} to hold. To this end, we first rearrange these conditions into a more convenient form.
We begin by combining the inequalities involving $1/r_1$ and $1/r_2$ in \eqref{c7}-\eqref{c9}. The first two conditions in \eqref{c8} reduce to
$$ {b+1} < \frac{n}{r_1} \leq \frac{(p-1)n}{r}$$
since $b>0$, $p>2$, and $r \ge 2$.
(The condition $r\ge2$ is guaranteed by the Schr\"odinger admissibility condition.)
Similarly, the first two conditions in \eqref{c9} reduce to 
$${b} < \frac{n}{r_2} \le \frac{pn}{r}.$$
Consequently, the conditions \eqref{c7}-\eqref{c9} reduce to 
\begin{equation}\label{c77}
\frac1{r_1}+\frac1{r_2}=\frac{1}{\tilde r'}+\frac{\alpha}{n}, \quad {b-p+2}=\frac{n}{r_1}-\frac{(p-1)n}{r}, \quad {b-p}=\frac{n}{r_2}-\frac{pn}{r},
\end{equation}
\begin{equation}\label{c88}
\frac{n}{r_1}<n,\quad {b+1}<\frac{n}{r_1}\leq \frac{(p-1)n}{r},\quad \frac{n}{r_2} < n, \quad {b}<\frac{n}{r_2}\le\frac{pn}{r}.
\end{equation}
We turn to the conditions \eqref{c10}-\eqref{c6}. 
Proceeding as before and using $b> 0$, $p>2$, and $r \ge 2$, we combine the inequalities involving $1/r_3$, $1/r_4$, $1/r_5$, and $1/r_6$ from \eqref{c11}-\eqref{c6} to obtain  
$${b} < \frac{n}{r_3} \leq \frac{(p-1)n}{r}, \quad 
{b+1} < \frac{n}{r_4} \leq \frac{pn}{r}, \quad 
b < \frac{n}{r_5} \leq \frac{(p-2)n}{r}, \quad 
{b}< \frac{n}{r_6} \leq \frac{(p-1)n}{r}.$$
Consequently, the conditions \eqref{c10}-\eqref{c6} reduce to
\begin{equation}\label{c100}
\frac{1}{r_3} + \frac{1}{r_4} = \frac{1}{\tilde r'} + \frac{\alpha}{n}, \quad {b-p+1}= \frac{n}{r_3} - \frac{(p-1)n}{r}, \quad {b-p+1}= \frac{n}{r_4} - \frac{pn}{r},
\end{equation}
\begin{equation*}
\frac{n}{r_3}<n, \quad {b} < \frac{n}{r_3} \leq \frac{(p-1)n}{r}, \quad \frac{n}{r_4}<n, \quad {b+1} < \frac{n}{r_4} \leq \frac{pn}{r},
\end{equation*}
\begin{equation*}
b < \frac{n}{r_5} \leq \frac{(p-2)n}{r}, \quad {b-p+2}=\frac{n}{r_5}-\frac{(p-2)n}{r},
\end{equation*}
\begin{equation}\label{c22}
b < \frac{n}{r_6} \leq \frac{(p-1)n}{r}, \quad {b-p+1}=\frac{n}{r_6}-\frac{(p-1)n}{r}.
\end{equation}

We now proceed to eliminate the auxiliary indices $r_1$, $\cdots$,$r_6$ from the reduced conditions \eqref{c77}-\eqref{c22} of  the conditions \eqref{c7}-\eqref{c6} required for \eqref{non11} to hold.
To eliminate $r_1$ and $r_2$, we substitute the last two conditions in \eqref{c77} into the first condition in \eqref{c77} and all the conditions in \eqref{c88}. This yields 
\begin{equation}\label{a}
\frac{(2p-1)n}{r}+2b-2p+2=\frac{n}{\tilde r'}+\alpha, \quad 1<\frac{n}{r} < 1+\frac{n-1-b}{p-1}, \quad b \le p-2,
\end{equation}
\begin{equation}\label{b}
\frac{n}r <1+\frac{n-b}{p} , \quad b\le p.
\end{equation}
Similarly, eliminating $r_3$, $r_4$, $r_5$, and $r_6$ from \eqref{c100}-\eqref{c22} yields
\begin{equation}\label{c16}
\frac{(2p-1)n}{r}+2b-2p+2=\frac{n}{\tilde r'}+\alpha, \quad 1<\frac{n}{r}<1+\frac{n-b}{p-1}, \quad b \le p-1,
\end{equation}
\begin{equation}\label{c}
\frac{n}{r}<1+\frac{n-1-b}{p}, \quad b\le p-2.
\end{equation}

Note that the last condition in \eqref{b} is directly implied by that in \eqref{a} and is therefore redundant. The last condition in \eqref{c16} is likewise implied by that in \eqref{c} and can be omitted.
Furthermore, the upper bound on $n/r$ in \eqref{c16} exceeds that in \eqref{a}, and is therefore not restrictive. Likewise, the upper bound on $n/r$ in \eqref{b} exceeds that in \eqref{c} and thus need not be considered.
Similarly, the upper bound on $n/r$ in \eqref{a} is greater than that in \eqref{c}, provided that $b<n-1$.
Since $p=2+\frac{\alpha-n+4-2b}{n-2}>2$, it follows that 
\begin{equation}\label{se}
b<\frac{4+\alpha-n}{2}<2,
\end{equation}
and hence $b<n-1$.
Consequently, after removing the redundant conditions from \eqref{a}-\eqref{c}, we are left with 
\begin{equation}\label{cd}
\frac{(2p-1)n}{r}+2b-2p+2=\frac{n}{\tilde r'}+\alpha, \quad 1<\frac{n}{r}<1+\frac{n-1-b}{p}, \quad b \le p-2.
\end{equation}

Next we take the intersection of the condition \eqref{cd} with \eqref{q} and the Schr\"odinger admissibility condition 
\begin{equation}\label{ad2}
0 \le \frac1q \le \frac12 , \quad \frac2q + \frac{n}r = \frac{n}2,\quad
0 \le \frac1{\tilde q} \le \frac12 ,\quad \frac2{\tilde q} + \frac{n}{\tilde r} = \frac{n}2.
\end{equation}
First we substitute \eqref{q} into the last two conditions of \eqref{ad2} to obtain
\begin{equation}\label{q1}
\frac1{2p-1} \le \frac{2}q \le \frac{2}{2p-1}, \quad 2-\frac{2(2p-1)}{q} +\frac{n}{\tilde r}=\frac{n}2.
\end{equation}
To eliminate $q$, we substitute the condition $2/q + n/r = n/2$ from \eqref{ad2} into the first condition in \eqref{ad2} and all the conditions in \eqref{q1}. This yields
\begin{equation}\label{r}
\frac{n-2}{2} \le \frac{n}{r} \le \frac{n}2, \quad \frac{n}{2}-\frac2{2p-1} \le \frac{n}{r} \le \frac{n}{2}-\frac1{2p-1}, \quad \frac{(2p-1)n}{r}+\frac{n}{\tilde r}= pn- 2. 
\end{equation}
Since $\frac{n-2}{2}<\frac{n}{2}-\frac2{2p-1}$ for $p>2$, the second condition implies the first one, and hence the first condition is redundant.
Substituting the last condition in \eqref{r} to the first condition in \eqref{cd}, it follows that $p=2+\frac{\alpha-n+4-2b}{n-2}.$ Hence these conditions can be removed.

Summarizing the outcome of the above reductions, the required conditions are
\begin{equation}\label{r0}
b\le p-2, \quad 1<\frac{n}{r}<1+\frac{n-1-b}{p}, \quad \frac{n}{2}-\frac2{2p-1} \le \frac{n}{r} \le \frac{n}2-\frac1{2p-1}.
\end{equation}
For such an $r$ to exist, the lower bounds on $n/r$ in \eqref{r0} must be  less than the upper bounds. This yields
\begin{equation*}\label{r1}
b<n-1, \quad p>\frac12 + \frac1{n-2}, \quad b < n-\frac{(n-2)p}{2}+ \frac{2p}{2p-1}-1,
\end{equation*}
which are all redundant. Indeed, we have shown above that the first condition $b<n-1$ is satisfied. The second condition is trivially satisfied since $p>2$ and $n\ge 3$.
Substituting $p=2+\frac{\alpha-n+4-2b}{n-2}$ into the last condition yields
$\frac{n-\alpha}{2}>\frac{n-2}{4b-2\alpha-2-n}$. 
Since $2b-\alpha<4-n\le1$ by \eqref{se}, the right-hand side is negative, and hence the condition is automatically satisfied.
Consequently, only the first condition in \eqref{r0} remains. Substituting $p=2+\frac{\alpha-n+4-2b}{n-2}$, we obtain
\begin{equation*}\label{r1'}
0< b\le \frac{\alpha-n+4}{n},
\end{equation*}
which coincides with the assumption \eqref{ass}.

Next, to prove \eqref{non22}, we use the following simple inequality:
\begin{align*}
\big|\mathcal{N}[u]-\mathcal{N}[v]\big| &= \Big||x|^{-b}|u|^{p-2}u(I_\alpha\ast|x|^{-b}|u|^p)-|x|^{-b}|v|^{p-2}v(I_\alpha\ast|x|^{-b}|v|^p)\Big|\\
&=\Big||x|^{-b}(|u|^{p-2}u-|v|^{p-2}v)(I_\alpha\ast|x|^{-b}|u|^p)\\
&\qquad\qquad\qquad\qquad\qquad\qquad\quad + |x|^{-b}|v|^{p-2}v\big(I_\alpha \ast|x|^{-b}(|u|^p-|v|^p)\big)\Big|\\
&\lesssim \Big||x|^{-b}(|u|^{p-2}+|v|^{p-2})|u-v||I_\alpha\ast|x|^{-b}|u|^p|\Big| \\
&\qquad\qquad\qquad\qquad\quad +\Big||x|^{-b}|v|^{p-1}\big(I_\alpha\ast|x|^{-b}(|u|^{p-1}+|v|^{p-1})|u-v|\big)\Big|.
\end{align*}
Consequently,
\begin{align*}
\big\| \mathcal N[u]-\mathcal N[v] \big\|_{L_t^{\tilde q'}(I;L_x^{\tilde r'})} &\lesssim \big\| |x|^{-b}(|u|^{p-2}+|v|^{p-2})|u-v|(I_\alpha\ast|\cdot|^{-b}|u|^p)\big\|_{L_t^{\tilde q'}(I;L_x^{\tilde r'})}\\
&\qquad \quad + \big\| |x|^{-b}|v|^{p-1}\big(I_\alpha \ast|\cdot|^{-b}(|u|^{p-1}+|v|^{p-1})|u-v|\big)\big\|_{L_t^{\tilde q'}(I;L_x^{\tilde r'})}.
\end{align*}
The first term on the right-hand side is estimated in the same way as $A_3$, replacing the factor $|u|^{p-2}\nabla u$ with $(|u|^{p-2}+|v|^{p-2})|u-v|$. The second term is treated analogously to $A_4$, replacing  $|u|^{p-1}\nabla u$ with $(|u|^{p-1}+|v|^{p-1})|u-v|$.
This yields
\begin{align*}
\big\| |x|^{-b}(|u|^{p-2}+|v|^{p-2})|u-v|&(I_\alpha\ast|\cdot|^{-b}|u|^p)\big\|_{L_t^{\tilde q'}(I;L_x^{\tilde r'})} \\
&\qquad \lesssim \big( \big\|\nabla u \big\|_{L_t^{q}(I;L_x^{r})}^{2p-2} + \big\| \nabla v \big\|_{L_t^{q}(I;L_x^{r})}^{2p-2}\big) \|u-v\|_{L_t^{q}(I;L_x^{r})}
\end{align*}
and 
\begin{align*}
\big\| |x|^{-b}|u|^{p-1}(I_\alpha\ast|\cdot|^{-b}(|u|^{p-1}+&|v|^{p-1})|u-v|)\big\|_{L_t^{\tilde q'}(I;L_x^{\tilde r'})} \\
&\quad \lesssim \big(\|\nabla u\|_{L_t^{q}(I;L_x^{r})}^{2p-2}+\|\nabla v\|_{L_t^{q}(I;L_x^{r})}^{2p-2}\big)\|u-v\|_{L_t^{q}(I;L_x^{r})}.
\end{align*}
Both estimates hold under the same conditions as those used in the treatment of $A_3$ and $A_4$, respectively.
Consequently, \eqref{non22} does not impose any additional condition.

\section{Well-posedness}\label{sec3}
Once the nonlinear estimates are established, the well-posedness result (Theorem \ref{thm}) follows in a standard way via the contraction mapping principle and the following Strichartz estimates.

\begin{lem}[\cite{St,gv,M,KT}]\label{str}
Let $n\geq3$. Assume that $(q,r)$ and $(\tilde q , \tilde r)$ are Schr\"odinger admissible pairs, that is,
\begin{equation*}\label{sch}
    2\le q,\, \tilde q \le \infty, \quad 
\frac{2}{q}+\frac{n}{r}=\frac{n}{2}, \quad
\frac{2}{\tilde q}+\frac{n}{\tilde r}=\frac{n}{2}.
\end{equation*} 
Then
\begin{equation}\label{homo}
\|e^{-it\Delta}f\|_{L_t^qL_x^r}\lesssim \|f\|_{L^2},
\end{equation}
\begin{equation*}\label{inhomo}
\Big\|\int_0^{t}e^{-i(t-s)\Delta}F(\cdot,s)ds \Big\|_{L_t^qL_x^r}\lesssim \|F\|_{L_t^{\tilde q'}L_x^{\tilde r'}}.
\end{equation*}
\end{lem}

We first rewrite the Cauchy problem \eqref{S} in integral form via Duhamel's formula:
\begin{equation}\label{du}
u(t)=e^{-it\Delta} u_0 + i\epsilon\int_0^t e^{-i(t-s)\Delta} \mathcal N[u](s,\cdot) ds:=\Phi(u),
\end{equation}
where $\mathcal N[u]=|x|^{-b}|u|^{p-2}(I_\alpha \ast |\cdot|^{-b}|u|^p)u$. 
For suitable values $T,M,N>0$, we shall show that $\Phi$ defines a contraction mapping on the space 
$$X(T,M,N)=\big\{u \in C_t(I;H_x^1) \cap L_t^{q}(I;H_x^{1,r}): \sup_{t\in I} \|u\|_{H_x^1}\leq M,\, \|u\|_{\mathcal{H} (I)}\leq N\big\},$$
equipped with the distance 
$$d(u,v)=\|u-v\|_{S(I)},$$
where $I=[0,T]$ and $(q,r)\in \mathcal{A}$.
We also define 
$$\|u\|_{\mathcal H(I)}:= \|u\|_{S(I)} + \| \nabla u\|_{S(I)}$$
and 
$$\|u\|_{\mathcal{H}'(I)}:= \|u\|_{S'(I)} + \| \nabla u\|_{S'(I)}.$$

First, we show that $\Phi$ is well defined on $X$.  
By Lemma \ref{str}, we get
\begin{equation}\label{w1}
\|\Phi(u)\|_{\mathcal{H}(I)}\lesssim \|e^{-it\Delta }u_0\|_{\mathcal{H}(I)} +\big\|\mathcal N[u]\big\|_{\mathcal{H}'(I)}.
\end{equation}
Moreover,
\begin{equation*}
\sup_{t\in I} \|\Phi(u)\|_{H_x^1}\leq C\|u_0\|_{H^1}+\sup_{t \in I}\left\|\int_0^t e^{-i(t-s)\Delta} \mathcal N[u](\cdot,s) ds\right\|_{H_x^1}.
\end{equation*}
To estimate the second term on the right-hand side of the above inequality, we proceed as follows.
Since $ \| u \|_{H_x^1}\lesssim \|u\|_{L_x^2} + \| u\|_{\dot H_x^1}$, and $e^{it\Delta}$ is a unitary operator on both $L^2$ and $\dot H^1$, we apply the dual estimate of \eqref{homo} in Lemma \ref{str} to obtain
$$\sup_{t\in I} \left\| \int_0^t e^{-i(t-s)\Delta} \mathcal N[u](\cdot,s) ds \right\|_{H_x^1} \lesssim \|\mathcal{N}[u]\|_{S'(I)} + \|\nabla \mathcal{N}[u]\|_{S'(I)}.$$
Hence, 
\begin{equation*}
\sup_{t\in I}\|\Phi(u)\|_{H_x^1} \lesssim \|u_0\|_{H^1}+\|\mathcal{N}[u]\|_{\mathcal{H}'(I)}.
\end{equation*}
On the other hand, by Proposition \ref{non}, we have 
\begin{align}\label{w2}
\nonumber
\big\|\mathcal{N}[u]\big\|_{\mathcal{H}'(I)}&\lesssim \big\|\nabla u\big\|_{S(I)}^{2p-1} + \big\|\nabla u\big\|_{S(I)}^{2p-2} \|u\|_{S(I)} \\
\nonumber
&= \big\| \nabla u \big\|_{S(I)}^{2p-2} \| u \|_{\mathcal{H}(I)}\\
&\leq C N^{2p-1}
\end{align}
for $u \in X$.
Substituting this into the previous estimate, we obtain
\begin{equation}\label{aa}  
    \sup_{t \in I}\|\Phi(u)\|_{H_x^1} \leq C \|u_0\|_{H^1} + CN^{2p-1}.
\end{equation} 
Next, for some $\varepsilon>0$ to be chosen later, we have
\begin{equation}\label{sm}
\|e^{it\Delta }u_0\|_{\mathcal{H}(I)}\leq \varepsilon,
\end{equation}
which holds for sufficiently small $T>0$ by the dominated convergence theorem.
Combining this with \eqref{w1} and \eqref{w2}, we conclude that
\begin{equation}\label{bb}
\|\Phi(u)\|_{\mathcal{H}(I)} \leq \varepsilon + CN^{2p-1}. 
\end{equation}
Finally, combining \eqref{aa} and \eqref{bb}, we see that $\Phi(u)\in X$ for $u \in X$ provided that
\begin{equation}\label{w3}
\varepsilon + CN^{2p-1} \leq N \quad \text{and} \quad C\|u_0\|_{H^1} + CN^{2p-1} \leq M.
\end{equation}

Next, we show that $\Phi$ is a contraction on $X$.
Using the same argument as in \eqref{w1}, we obtain
\begin{equation*}
\big\| \Phi(u)-\Phi(v) \big\|_{S(I)} \lesssim \big\| \mathcal{N}[u]-\mathcal{N}[v] \big\|_{S'(I)}.
\end{equation*}
Applying Proposition \ref{non} (see \eqref{non2}), we have
\begin{align*}
\big\|\mathcal{N}[u]-\mathcal{N}[v] \big\|_{S'(I)}&\leq C \Big(\big\|\nabla u\big\|_{S(I)}^{2p-1}+\big\| \nabla v \big\|_{S(I)}^{2p-2}\Big)\|u-v\|_{S(I)} \\
&\leq C N^{2p-2}\|u-v\|_{S(I)}
\end{align*}
as in \eqref{w2}.
Therefore, for $u,v \in X$, we obtain 
$$d\big(\Phi(u), \Phi(v)\big)\leq C N^{2p-2}d(u,v).$$
Now, taking $M=2C\|u_0\|_{H^1}$ and $N=2\varepsilon$, and choosing $\varepsilon>0$ sufficiently small so that  \eqref{w3} holds and $CN^{2p-2}\leq 1/2$, it follows that $\Phi$ is a contraction on $X$.
Hence, by the contraction mapping principle, there exists a unique local solution $u \in C_t(I;H_x^1) \cap L_t^q(I;H_x^{1,r})$ for any $(q,r)\in \mathcal{A}$.

Using the homogeneous estimate \eqref{homo} in Lemma \ref{str}, we also observe that \eqref{sm} is satisfied if $\|u_0\|_{H^1}$ is sufficiently small:
\begin{equation*}
\|e^{-it{\Delta}}u_0\|_{\mathcal{H}(I)} \leq C \|u_0\|_{H^1}\leq \varepsilon.
\end{equation*}
In this case, one can take $T=\infty$ in the above argument, which yields a global unique solution. 

The continuous dependence of the solution $u$ on the initial data $u_0$ follows in the same way. Indeed, we have
\begin{align*}
d(u,v) & \lesssim d (e^{-it\Delta} u_0 , e^{-it\Delta} v_0 ) + d\Big(\int_0^t e^{-i(t-s)\Delta}\mathcal{N}[u]ds , \int_0^t e^{-i(t-s)\Delta}\mathcal{N}[v]ds \Big) \\
&\lesssim \|u_0 - v_0\|_{L^2} + \frac12 d(u,v).
\end{align*}
Rearranging the above inequality, we obtain
\begin{align*}
d(u,v) \lesssim \|u_0-v_0\|_{H^1}.
\end{align*}
Here, $u$ and $v$ denote the corresponding solutions associated with the initial data $u_0$ and $v_0$, respectively.

To prove the scattering property, using \eqref{du} and \eqref{w2}, we first note that 
\begin{align*}
\big\| e^{it_2 \Delta}u(t_2)-e^{it_1 \Delta}u(t_1) \big\|_{H_x^1} = \Big\| \int_{t_1}^{t_2} e^{is\Delta}\mathcal{N}[u]ds\Big\|_{H_x^1}
& \lesssim \big\| \mathcal{N}[u] \big\|_{\mathcal{H}'([t_1 ,t_2])}\\
& \lesssim \|u\|^{2p-1}_{\mathcal{H}([t_1,t_2])} \rightarrow 0
\end{align*}
as $t_1, t_2 \rightarrow {\infty}$.
This implies that $\phi :=\lim_{t\rightarrow {\infty}} e^{it \Delta} u(t)$ exists in $H^1$.
Using this and \eqref{du}, we write
$$u(t)-e^{-it\Delta}\phi=-i\varepsilon\int_t^{\infty} e^{-i(t-s)\Delta}\mathcal{N}[u]ds,$$
and hence 
\begin{align*}
\big\|u(t)-e^{-it\Delta}\phi \big\|_{H_x^1}=\Big\|\int_t^{\infty} e^{-i(t-s)\Delta}\mathcal{N}[u] ds\Big\|_{H_x^1} 
&\lesssim\|\mathcal{N}[u]\|_{\mathcal{H}'([t,\infty])} \\
&\lesssim \|u\|^{2p-1}_{\mathcal{H}([t,\infty])} \rightarrow 0
\end{align*}
as $t\rightarrow {\infty}$. 

This completes the proof of the theorem.


\end{document}